\documentclass[letterpaper, 10 pt, conference]{ieeeconf}  

\IEEEoverridecommandlockouts                              

\usepackage{graphics} 
\usepackage{epsfig} 
\usepackage{mathptmx} 
\usepackage{times} 
\usepackage{amsmath} 
\usepackage{amssymb}  
\DeclareMathOperator{\E}{\mathbb{E}}
\newcounter{cor}
\newtheorem{cor}{Corollary}
\newcounter{remark}
\newtheorem{remark}{Remark}
\newcounter{example}
\newtheorem{example}{Example}
\newcounter{definition}
\newtheorem{definition}{Definition}
\newcounter{lemma}
\newtheorem{lem}[lemma]{Lemma}
\newcounter{thm}
\newtheorem{thm}{Theorem}

\title{\LARGE\bf
Optimal Contract Design for Incentive-Based
Demand Response}

\author{Donya G. Dobakhshari and Vijay Gupta \thanks{The authors are with the Department of Electrical Engineering, University of Notre Dame, IN, USA. Email:  (\texttt{dghavide, vgupta2)@nd.edu.} The work was supported in part by NSF award 1550016, 1239224 and 1544724. This work has been submitted to the IEEE for possible publication. Copyright may be transferred without notice, after which this version may no longer be accessible.}}

\begin{document}

\maketitle
\thispagestyle{empty}
\pagestyle{empty}

\begin{abstract}
We design an optimal contract between a demand response aggregator (DRA) and a customer for incentive-based demand response. We consider a setting in which the customer is asked to reduce her consumption by the DRA and she is compensated for this reduction. However, since the DRA must supply the customer with as much power as she desires, a strategic customer can temporarily increase her base load to report a larger reduction as a part of the demand response event. 
The DRA wishes to incentivize the customer  both to  make the maximal effort in reducing the load and to not falsify the base load. We model the problem of designing the contract by the DRA for the customer as a management contract design problem and present a solution. The optimal contract consists of two parts: a part that depends on (the possibly inflated) load reduction as measured and another that provides a share of the  profit that ensues to the DRA through the demand response event to the customer.
\end{abstract}

\section{INTRODUCTION}
\label{i0}
Demand response, in which a utility company or an aggregator motivates customers to curtail their power usage, has now become an acceptable method in situations when high peaks in demand occur.  Demand Response (DR) can be defined as the change in electric usage by end-use customers from their normal consumption patterns in response to changes in the price of electricity or any other incentive~\cite{va, deng20} and~\cite{imp6}.

Generally, DR programs are divided into two main categories: Incentive Based Programs (IBP) and Price Based Programs (PBP).  PBPs provide time of usage based electricity prices and the consumers are expected to adjust their demand in response to such a price profile. On the other hand, IBPs offer  incentives to customers to reduce their demand. These incentives may be constant and based only on customer participation in the program (classical) or  dynamic in the sense that they vary with the amount of load reduction that a customer achieves (market-based). There exists a rich literature (e.g, in~\cite{qu,mo,sam,ro, wa} and the references therein) studying issues such as social welfare maximization, minimization of electricity generation and delivery costs, and reducing renewable energy supply uncertainty for incentive-based demand response.

In this paper, we consider an incentive based DR scenario where  participants are rewarded financially by the demand response aggregator (which role can also be filled by a utility company) for the amount of load reduction provided by consumers during DR events. However, unlike the existing literature, we consider a strategic customer that maximizes her own profit by predicting the impact of her actions and the information she transmits. Specifically, by taking advantage of the fact that the demand response aggregator (DRA) must supply as much power as the customer desires, a strategic customer can artificially inflate her base load before an expected DR event. Then, during the DR event, for the same \textit{nominal} load reduction,  the customer can report more \textit{measured} load reduction and gain more financial reward from the DRA.

In such a scenario, we wish to find a contract that incentivizes the strategic customer to achieve the maximum nominal load reduction possible. The main contribution of our work is to characterize an optimal contract for this DR problem. Our solution is similar to a managerial contract  model studied e.g. in ~\cite{imp3, imp2}; however, we do not assume accurate knowledge of the profit achieved by the DRA  as a result of the load reduction by the customer.  The optimal contract consists of two parts: a part that depends on the reported load reduction and another that provides a share of the  profit for the DRA through the demand response event to the customer.

One interesting result is that the optimal contract leads to \textit{under-reporting} of load reduction by the customer up to a specific value of nominal load reduction and  \textit{over-reporting} of reduction above that value. In other words, if the strategic customer wishes to maximize her profit, she may sometimes decrease her base load before the DR event to under-report her power reduction.  Furthermore, if the expected difference between the nominal reduction with the true base load and the reported one with the inflated or deflated base load is positive, the DRA's expected profit is an increasing function with respect to the  share provided  to the customer.  Finally, the analysis implies that it is always optimal (from the DRA's viewpoint) to assign some positive  share of the profit to the customer.

The paper is organized as follows. In Section \ref{i1}, the problem statement is presented. In Section \ref{i2}, the solution to the optimization problem is presented. Next, we discuss the optimal  contract structure and its properties in Section \ref{i01}.  The final section concludes the paper by pointing out some directions for future work.

\paragraph*{Notation} $f_{X|A}(x|a)$  (which is often simplified to $f(x|a)$) and $F_{X|A}(x|a)$ denote the probability distribution function (pdf) and cumulative distribution function (cdf) of random variable $X$ given the event $A=a$ respectively. Gaussian distribution is denoted by  $\mathcal N(m,\sigma^2 )$ where $m$ is the mean and $\sigma $ is the standard deviation. Derivative of a function $W$ with respect to a variable $x$ is denoted as $W_x$ or $W'$ if the variable is clear from the context. For two functions $g$ and $h$, $g*h$ denotes the convolution between $g$ and $h$. $\E[Y]$ denotes expectation of random variable $Y$. By abusing notation, we sometimes write the expectation as $\E[y]$.

\section{Problem Statement }
\label{i1}
We model a demand response event as beginning when the DRA calls on a customer to decrease her power consumption. A strategic customer, anticipating such a call, can increase her \textit{base load}, or the load before the demand response event. This pre-increase allows the customer to reduce the load by a larger amount than would have been possible in the absence of such an increase. The customer potentially gains from this larger reduction if the market based DR entails payment of an incentive proportional to the load reduction by the customer during the DR event.  On the other hand, a contract must make the payment proportional to the load reduction to exert the maximal effort for reducing the load by as much amount as possible (\textit{See Example 1 below}). 

\textbf{Remark:} It is worth pointing out that the falsification of the load reduction reported to the DRA happens even though the load at the customer is being monitored constantly and accurately. Further, the DRA can not find the `true' base load by considering the load used by the customer at some arbitrary time before the DR event. For one, this simply shifts the problem of customer manipulation of the load to an earlier time. Second, some of the increase in the base load may be due to true shifts in customer need due to, e.g., increased temperature.

\subsection{Problem Formulation}

\begin{figure}[tbp]
\noindent%
\includegraphics[width=0.9\linewidth]{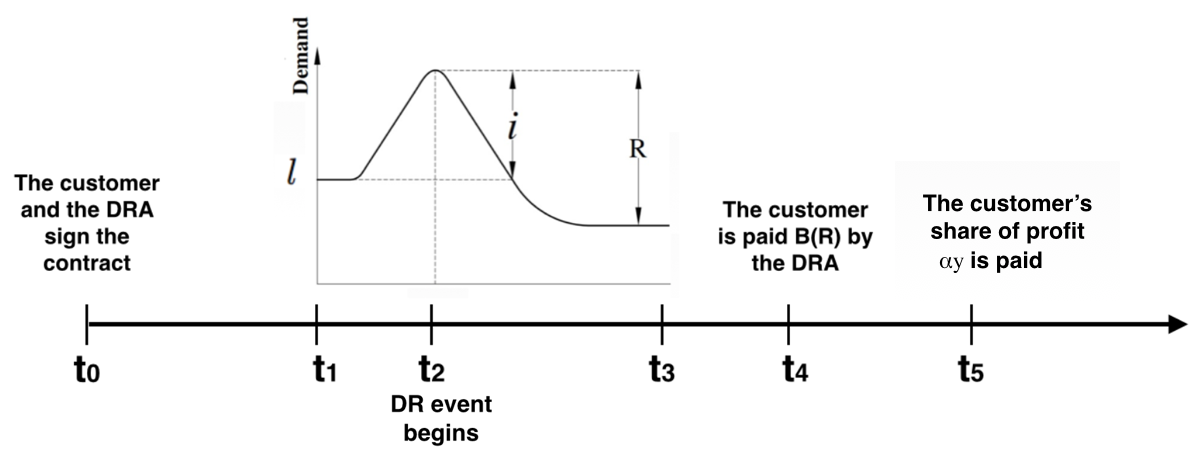}
\caption{Timeline of the DR event and the proposed contract.}
\label{pic_2}
\vspace{-1mm}
\end{figure}

Refer to the timeline shown in Figure \ref{pic_2}. The true base load (without any manipulation) is given by $l$. At time $t_1$, the customer calculates the effort  $a$ she is willing to put in for achieving the load reduction $x$ during the DR event. The load reduction is according to the probability density function  $f(x|a)$ which is public knowledge. We assume that an effort $a$ costs the customer $h(a)$ ($h(a)$ is convex and $h(0)=0$). Further, this effort and the planned load reduction $(a,x)$ depends on private knowledge at the customer and hence can be calculated by the customer, but not by the DRA. For instance, a factory might be able to induce a large load reduction with a small effort based on its assembly line requirements given the orders it has to fulfill. After this calculation, the customer at time $t_{1}$ may increase (or decrease) the load by an amount $i$ in anticipation of the DR event. 

At time $t_2$, the DR event begins and the DRA calls on the customer to decrease her load. The customer now makes the effort $a$ yielding a reduction of the load by $x$. The DR  event ends at $t_3$ with the customer having decreased the load by an amount $R(x)$ (which is often simplified to $R$). Note that the planned reduction in the load is  $x= R(x)-i(x)$, while the false reported load reduction  
 is $R(x)$. We also show the times $t_0$, $t_4$ and $t_5$ in the timeline in Figure \ref{pic_2}. At time $t_0$ (much before $t_{1}$), the contract is signed between the DRA and the customer, while at times $t_4$ and $t_5$, the customer is paid by the DRA according to the contract we will propose in the sequel. We note that $t_{0}$ is sufficiently early, so that at $t_{0}$, the customer does not know the local conditions and must consider her expected utility according to the probability density function $f(x|a)$. 

The time $t_{4}$ is sufficiently close to the DR event, so that, the realized value of $x$ is not known at time  $t_4$ to the DRA. The customer needs to be paid at least in part at $t_{4}$ to incentivize her to participate in the DR event. However, at some (much) later time  $t_5$, the DRA may be able to estimate the realized value of $x$, possibly with some error. This noisy estimate can be obtained by, e.g., large scale data analysis on all similar customers on that day or historical behavioral of the same customer. We denote this estimate by $y$ where $y=x+n$. We assume that the random variables $X$ and $N$ are independent and $N\sim \mathcal N(m_n,\sigma_n ^2)$.

We model  the falsification cost incurred by the customer (e.g. extra charge paid for boosting her consumption) by a quadratic function for simplicity and denote this cost as $g(R-x)=\frac{(R-x)^2}{2}$. 

The problem is for the DRA to design a contract  that maximizes his own profit. Since this profit depends on the load reduction by the customer, the contract must induce a rational customer to choose an \textit{action} $a$ and a \textit{report} $R$ that are optimal for the DRA. The profit for the DRA occurs due to the load reduction by the customer, modulo any payments made to the customer as part of her contract. We discuss the intuition for the proposed contract through some examples.  
\begin{example}
Consider a contract that provides a constant incentive $c$ to the customer for decreasing her load. Then, the customer's utility is given by:
\begin{displaymath}
V_{cust} = cu(R)-g(R-x)-h(a),
\end{displaymath}
where $u(.)$ is the unit step function. The DRA's utility 
is given by
\begin{displaymath}
\Pi_{DRA}=y-cu(R).
\end{displaymath}
In this case, the customer seeking to maximize her utility, independent of the value of $c$, will choose $a=0$ and respectively $x=0$ (i.e., no action and no true load reduction) but $R=0^{+}$ (i.e., minimal load reduction). This implies that to induce positive load reduction, a contract must make at least part of the payment proportional to the load reduction.
\end{example}
\begin{example}
Consider a contract in which the DRA provides an incentive $cR$ to the customer in response to the reported reduction $R$ (which is all that she has access to at $t_{4}$). Then, the customer's utility is given by:
\begin{displaymath}
V_{cust} = cR-g(R-x)-h(a),
\end{displaymath}
while the DRA's utility 
is given by
\begin{displaymath}
\Pi_{DRA}=y-cR.
\end{displaymath}
The DRA seeks to optimize $\Pi_{DRA}$ over $c$ assuming that the customer will choose $a$ and $R$ to maximize $V_{cust}$. However, irrespective of the optimization, a  customer again can take no action, i.e. $a=x=0$, and report $R=c$ to gain the positive profit ${c^2}/{2}$. Thus, a good contract  must entail some payments that depends on the DRA's estimated  profit $y$.
\end{example}

Next, we propose a contract structure free from the shortcomings of these intuitive contracts.


\subsection{Contract Structure}

Inspired by managerial contracts studied e.g. in \cite{imp3, imp2}, we propose a contract of the form $(\alpha, B(R))$ in which $\alpha$ refers to the share of his own profit that the DRA  provides to the customer, while $B(R)$ refers to the payment made in proportion to the reported reduction $R(x)$. Referring to Figure \ref{pic_2}, to incentivize the customers to participate in the program, $B(R)$ is paid at $t_{4}$. However, the shares (even though they are allotted at $t_{0}$) can be encashed only at a much later time $t_{5}$ when an estimated value of the profit can be calculated and revealed. Note that the portion of the payment corresponding to the share $\alpha$ is calculated on the basis of the noisy estimate $y$ of $x$. 
Thus, the customer's utility is given by 
\begin{equation}
V\stackrel{\triangle}{=}  \alpha y-g(R(x)-x)- h(a)+B(R(x)),
\end{equation}
 while the DRA utility is given by
 \begin{equation}
 \Pi \stackrel{\triangle}{=} (1-\alpha)y- B(R(x)).
 \end{equation} It is worth pointing out that  the customer will report a load reduction to realize  $x$; therefore, $R$ is a function of $x$, not $y$.

 Thus, the optimization problem to be solved by DRA (subject to participation and rationality constraints for customer) is given by\begin{equation}
\underset {\alpha , B(R)} {\max } \E[\Pi]=\underset {\alpha , B(R)} {\max } \E[(1-\alpha)y- B(R(x))]
\label{e1}
\end{equation}

We now describe the constraints for the optimization problem in (\ref{e1}). 

\begin{enumerate}
\item Rationality in the choice of effort: The first assumption is that the customer chooses the level of effort $a$ to maximize her expected utility $\E[V]$. Thus, the first two constraints are given by  $\frac{\partial \E[V]}{\partial a}=0 $ and $\frac{\partial^2\E[V]}{\partial a^2}\leq0$ where the expectation is taken with respect to $x$ and $y$.
\item Ex ante individual rationality: The expected utility of the customer must be positive to ensure that she participates in the DR event. This implies a constraint of the form $\E[V]\geq0$.
 \item Interim individual rationality: 
 We impose that the customer must be incentivized to continue even though she can choose to leave after she makes effort $a$ and sees her comfort reduced. We impose this constraint as 
 \begin {equation} W\stackrel{\triangle}{=} V+h(a)\geq0.
  \end {equation} 
 
 \item Incentive compatibility: We impose two further constraints $R'(x)\geq0$ and $W_x=\alpha+g'(R(x)-x))$ as incentive compatibility constraints that ensure truthtelling by the customer in the conditional direct revelation contract~\cite{imp3}. 
\end{enumerate}
Thus, the optimization problem can be rewritten as 

\begin{equation}
\underset {\alpha , B(R)} {\max } \E[\Pi]=\underset {\alpha , B(R)} {\max } \E[(1-\alpha)y-B(R)]=\underset {\alpha , B(R)} {\max }\E[y-g-W]
\label{e3}
\end{equation}

subject to:
\begin{subequations}
\label{don}
\begin{align}
 \frac{\partial \E[W]}{\partial a}-h'(a)=0, \hspace{2mm} \frac{\partial^2 \E[W]}{\partial a^2}-h''(a)\leq0 \label{e03} \\ 
  \E[W]-h(a)\geq0 \label{e0}  \\ 
 W\geq0 \label{e04} \\
  W_x=\alpha+g'(R(x)-x) , \hspace{2mm} R'(x)\geq0  \label{e05}
\end{align}
\end{subequations}

We will present the optimal contract in Section \ref{i2}. The solution depends on the properties of the pdfs $f_{X|A}(x|a)$ that describes the planned reduction $x$ based on effort $a$ of the customer and $f_N(n)$ which is the pdf of the estimation error in the knowledge of $x$ in the long term. We make the following assumptions about these functions :

\begin{enumerate}
\item Assumption $\mathcal{A}_{1}$:  The cdf $F_{X|A}(x| a)$ of $f_{X|A}(x|a)$ is strictly decreasing, convex and continuously differentiable in $a$ for all $x$ and for all $a$. This is a natural assumption implying that higher customer effort induces a first-order stochastic improvement in the distribution of load reduction and results in diminishing marginal returns from effort.
\item  Assumption $\mathcal{A}_{2}$: $f_{X|A}(x|a)>0$ for all $x$ and $a$. Further, there exists $M>0$ such that for all $x$ and $a,$ $f_{X|A}(x|a)<M$ and $\frac{\partial{f_{X|A}(x|a)}} {\partial{x}}< M$. 
\item Assumption $\mathcal{A}_{3}$: $\frac{\partial{f_{X|A}(x|a)}} {\partial{x}}\geq0$ for all $a>0$.  This assumption implies that positive values of load reduction are more likely than zero values of load reduction.
\item Assumption $\mathcal{A}_{4}$: $(F_{X|A}(x|a)-1)/ f_{X|A}(x|a)$ is strictly concave in $x$ for all $a$ and $F_a(x|a)/ f_{X|A}(x|a)$ is strictly convex in $x$ for all $a$ where  
\begin{equation}
F_a(x|a)=\frac{\partial{F_{X|A}(x|a)}} {\partial{a}}.
\end{equation}
 Assumption $\mathcal{A}_{4}$ is provided for the proof of  corollary \ref{pou} (in section \ref{i01}) where we restrict to Gaussian distribution for $x$ and $n$ for simplicity.
 \end{enumerate}
 
 \begin{lem}
 The properties encapsulated in assumptions $\mathcal{A}_{1}$-$\mathcal{A}_{4}$ hold for the probability density function $f_{Y|A}(y|a)$ and the corresponding CDF $F_{Y|A}(y|a)$ as well.
 \end{lem}

 \begin{proof}
We present the proof for assumption $\mathcal{A}_{1}$, the proofs for rest of assumptions are similar.

By definition, $x$ and $y$  are related as $y=x+n$. Thus, $F_{Y|A}(y)$ and $f_{Y|A}(y)$ can be derived as follows:
\begin{equation}
f_{Y|A}(y|a)= f_{X|A}(y|a)*f_N(y)
\end{equation}
\[
F_{Y|A}(y|a)= F_{X|A}(y|a)*f_N(y)
\]
\begin{equation}
{\partial F_{Y|A}(y|a)}/{\partial a}= [{\partial F_{X|A}(y|a)}/{\partial a}]*f_N(y).
 \label{oo}
\end{equation}
According to \eqref{oo} and noting that $f_N(n)$ is a probability distribution function and positive everywhere,
\begin{equation}
\frac{\partial(F_{X|A}(y|a))}{\partial a}\leq0 \rightarrow \frac{\partial(F_{Y|A}(y|a))}{\partial a}\leq 0. 
\end{equation}
Thus, if $F_{X|A}(x|a)$ is strictly decreasing, $F_{Y|A}(y|a)$ will be strictly decreasing too. For convexity, since $F_{Y|A}(y|a)$ is strictly decreasing, it is enough to only prove  $\frac{\partial^2(F_Y(y|a))}{\partial a^2}\geq 0$. Now, given  \eqref{oo}, 
\begin{equation}
\frac{\partial^2(F_{X|A}(y|a))}{\partial a^2}\geq0   \rightarrow\frac{\partial^2(F_{Y|A}(y|a))}{\partial a^2}\geq 0.
\end {equation}
Therefore, Assumption $\mathcal{A}_{1}$ holds for $F_{Y|A}(y|a)$.\\

 \end{proof}


 \section{Optimal contract Structure }
 \label{i2}
 In this section, we present the solution of the optimization problem stated in~\eqref{e3}. 
Consider the argument being optimized in~\eqref{e3}. We begin with the case when $m_n=0$. We can rewrite the expected utility  $\Pi$ of the DRA as:
\begin{equation}
\E[y-g-W]=\E[x-g-W]=\int_{}^{}(x-g-W)f(x|a) dx.
\end{equation}
 We can define $U(W,R,x)=(x-g-W)f(x|a)$ to rewrite the optimization problem as\[
\underset {\alpha , B(R)} {\max } \Pi=\underset {\alpha , B(R)} {\max}\:\int_{}^{} U(W(x),R(x),x) dx,
\]
 subject to $W_x=\alpha+g'(R(x)-x),$ and the constraints in~\eqref{e03}-\eqref{e05}. This is an optimal control problem where the state variable is $W$ and $R$ is the control input.  We can solve this problem using the standard Hamiltonian approach. For now, we drop the second constraint \eqref{e05}  and will add this constraint later to the contract. Also we note the following result that was  proved in~\cite{imp3} and implies that the second order condition in \eqref{e03} is non-binding.
 \begin{lem}
 Given the distribution assumptions $\mathcal{A}_{1}-\mathcal{A}_{3}$, $\frac{\partial^2(E(W))}{\partial a^2}-h''(a)$ is strictly negative for any optimal contract.\end{lem} 
 
 Thus, we can  form the Hamiltonian 
\begin{equation}
 H=(y-g(R-x)-W)f+\varphi(\alpha+g'(R-x)), 
 \label{t2}
 \end{equation}
 and the corresponding Lagrangian 
\[
L=H+\tau W-\mu[\alpha+g'(R-x))F_a+h'(a)f]+\lambda (W-h)f,
\]
where $\varphi$ is the co-state variable (considering $W_x=\alpha+g'(R(x)-x)$), $f$ is $f(x|a)$, $W$ is the state variable, and $R$ is the control input. Further, $\mu$, $\lambda$, and $\tau$ are all non-negative multipliers included for considering the constraints \eqref{e03}, \eqref{e0}, and  \eqref{e04} respectively. 

We can now provide the structure of the optimal contract in the following result.
\begin{thm}
Given the assumptions $\mathcal{A}_{1}$-$\mathcal{A}_{3}$, if there exist a piecewise continuous function $\varphi(x)$, constraint multipliers $\mu$, $\lambda$ and $\tau$, and a  contract $(\alpha, B(R))$ that satisfy
: 

\begin{subequations}
\begin{align}
R-x=(\varphi-\mu F_a(x|a))/f(x|a) \label{b1}\\
-\varphi'=-(1-\lambda)f(x|a)+\tau \label{55}\\
a\left [\int_{0}^{1}(\alpha+g'F_a(t|a) dt+h'(a)\right]=0\\
\lambda(E(W)-h(a))=0\:  \textrm{ and }\lambda \geq0 \\
\quad \varphi(0)\leq0 \:, \varphi(1)\geq0, \hspace{2mm} \varphi(0)W(0)=\varphi(1)W(1)=0\\
\tau W=0  \textrm{ and }\: \tau \geq0 \\
a= \arg\max_a \E[y-g(R(x)-x)-W]\label{b3}, \hspace{2mm} R'(x)\geq 0,
\end{align}
\end{subequations}
then $(\alpha, R)$ is an optimal conditional contract that solves the optimization problem in \eqref{e3} and \eqref{don}.  \end{thm}
\begin{proof}
The proof follows readily from \cite[Chapter 6, Theorem 1]{notes} by substituting $W\rightarrow x $, $x\rightarrow t $, $R\rightarrow u$.
\end{proof}

This characterization can be used to determine the optimal values of the share $\alpha$ and bonus $B(R )$ on one hand, and the resulting effort $a$ and reporting function $R$ that are induced on the other. For illustration, we present the result for the reporting function below. The results for the other quantities can be derived similarly; we present insights on their forms in the next section.
\begin{definition}
Define $\hat{x}$ as the solution of the equation $W_{x}=0$. 
\end{definition}
The first equation in \eqref{e05} and the fact that the function $g(.)$ is a quadratic function implies that $\hat{x}$ can be obtained as
\begin{equation}
\hat{x}=R+\alpha. \label{me}
\end{equation}
 The reporting function induced by the optimal contract  is presented in the following result. 
\begin{thm}
The optimal reporting function is given as   \begin{equation}
R(x)=\begin{cases}
x-\alpha&\text{if} \  x\leq\hat{x}\\
x+\frac{(1-\lambda)(F_{Y|A}(x|a)*f_{N}(x)-1)-\mu {F}_a(x|a)*f_{N}(x)}{f_{Y|A}(x|a)*f_{N}(x)} & \text{if} \  x>\hat{x}\\
\end{cases}
\label{e6}
\end{equation} 
\end{thm}
\begin{proof}
The proof follows along the lines outlined in~\cite{imp3} using the conditions in \eqref{b1}-\eqref{b3} and the fact that $f_{N}(x)=f_{N}(-x)$.
\end{proof}

\begin{figure}[tbp]
\centering
\includegraphics[width=0.675\linewidth]{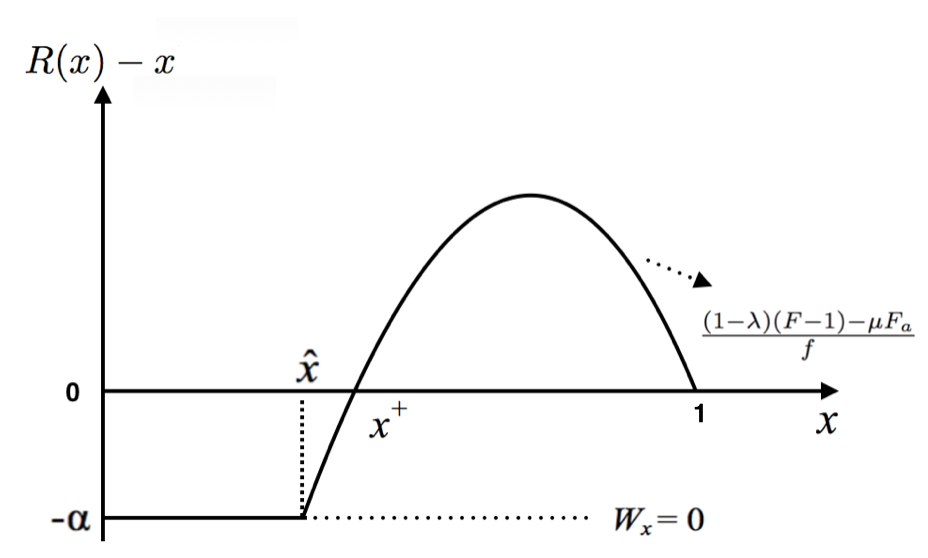}
\caption{The form of the reporting function.  There exists an $x^+$ such that  $R-x>0$ for $x>x^+$ (i.e., over-reporting happens) and $R-x<0$ for $x<x^+$ (i.e., under-reporting happens). The falsification is $-\alpha$  for $x\leq\hat{x}$ and $\frac{(1-\lambda)(F-1)-\mu F_a}{f}$ for $x>\hat{x}$. }
\label{pic_1}
\vspace{-2mm}
\end{figure}

\section{Discussion of Results}
\label{i01}
We now interpret the  results  obtained  in the previous section. For ease of interpretation and without loss of generality, we scale $x$ down to the range $[0,1]$.
\subsection{Form of the reporting function}
We can obtain a clearer interpretation of over-reporting (inflation of base load) and underreporting (reduction of base load) through the following result that specifies the form of the reporting function.
\begin{cor}
There exists $x^{+}>\hat{x}$, $x^{+}<1$ with $\hat{x}$ given in \eqref{me} such that the  optimal reporting function satisfies the following relation 
\begin{equation}
\begin{cases}
R<x \:&\text{if}\ x\leq\hat{x}\\
R<x \:&\text{if}\ \hat{x}<x\leq x^+\\
R>x \: &\text{if} \ x>x^+
\label{ss}
\end{cases}
\end{equation}
\end{cor}
\begin{proof}
From equation \eqref{e6}, we see that for $x\leq\hat{x}$, $R=x-\alpha$, which is always less than $x$. As $x$ increases, we appeal to assumption $\mathcal{A}_{4}$, its generalization to $y$ in Lemma 1 and  the fact that $f_{N}(x)=f_{N}(-x)$, to obtain that the function  \begin{equation*}
\frac{(1-\lambda)(F_{Y|A}(x|a)*f_{N}(x)-1)-\mu F_{{Y|A}_a}(x|a)*f_{N}(x)}{f_{Y|A}(x|a)*f_{N}(x)}
\end{equation*} is strictly concave in $x$ for all $a$. Thus, $R-x$ is an increasing function of $x$ and for a high enough value of $x$, the sign of $R-x$ will become positive~\cite{imp3}. This value is $x^{+}$ which is clearly larger than $\hat{x}$. 
 \end{proof}

The form of the reporting function is illustrated in Figure \ref{pic_1}. The result clarifies how the customer will falsify the load reduction by changing the base load. For nominal load reduction above $x^+$, $R-x>0$, i.e. the customer first increases the base load and then lowers it by the amount $R$. However, if $x<x^+$, $R-x<0$. This implies that in this case, the customer \textit{lowers} the demand at the beginning (or reports that she was going to reduce the demand even without the DR event) and then decreases the demand by $R$ again when called. This non-intuitive behavior can be understood if we remember that although $B(R)$ granted to customer is decreased through under-reporting, the share $\alpha$ of the profit assigned  to the customer to incentivize her to participate is larger in this case and this share compensates for the decrease in $B( R)$.


\subsection{Optimal Compensation}
In order to study the optimal compensation, we first present the following result without proof.

\begin{lem}In an optimal contract,  $W$  and $W_x$ are 0 for $x\leq \hat{x}$   and $W$ is greater than zero for $x>\hat{x}$.
 \end{lem}
We notice that the bonus can be considered to be a function of the savings $x$ and written as $B(x)$. Further, the bonus is related to $W$ as 
\begin{equation}
B(x)=W(x)+g(R(x)-x)-\alpha y.
\end{equation}
\begin{cor}
The optimal $B(x)$ satisfies the relation
\begin{equation}
\begin{cases}
B'(x)\leq0 &\text{if}\ x\leq\hat{x}\\
B'(x)\leq0 &\text{if}\ \hat{x}<x\leq x^+\\
B'(x)\geq0 &\text{if} \ x<x^+.
\end{cases}
\end{equation}
\end{cor}
\begin{proof}
When $x\leq \hat{x}$, $W(x)=0$. Since, $B(x)=g(-\alpha)-\alpha (x+n)$, we observe that $B'(x)=-\alpha$. Similarly, for $x> \hat{x}$
  \begin{equation}
  B(x)=g(R(x)-x)+\int_{\hat{x}}^{x} g'(R(t)-t)dt -\alpha \hat{x}-\alpha n.
  \label{e23}
  \end{equation} 
Using this result along with  the relation $R-\hat{x}=-\alpha$ leads to
\begin{equation}
 B'(x)=g'(R-x)R'(x)=(R-x)R'(x).
 \end{equation} Combining the two, cases we have 
\begin{equation}             
B'(x)=\begin{cases}
-\alpha&\text{if}\ x\leq\hat{x}\\
 (R-x)R'(x)&\text{if} \ x>\hat{x}.
\label{s4}
 \end{cases}
  \end{equation}
While $B'(x)<0$ if $x\leq\hat{x}$, for $x>\hat{x}$, the sign of $B'(x)$ depends on the sign of $R-x$  (since $R'(x)\geq0$ for $x>\hat{x}$). Thus, combining \eqref{s4} with \eqref{ss}   yields the desired result.
   \end{proof}
This result once again sheds light on the structure of  the two counteracting incentives provided to the customer. As $x$ increases, the bonus decreases up to the level $x^{+}$. In this range, the customer chooses to rely on the long term share and under-reports the load reduction she has made. For $x$ large enough, the bonus is an increasing function. In this range, the bonus is large enough and hence the customer chooses to boost her bonus by over-reporting her load reduction.

\subsection{Impact of Estimation Error}
In order  to compare the optimal reporting as a function of the noise in the estimation of the profit made due to the reduction of load, we need to investigate the optimal reporting function for the cases when $x$ is realized at $t_{5}$ exactly and with some error. Equation \eqref{e6} shows the relation between the optimal reporting function and the true profit with estimation error. In the absence of any error, the expression reduces to\begin{equation}
R(x)=\begin{cases}
x-\alpha&\text{if} \  x\leq\hat{x}\\
x+\frac{(1-\lambda)(F(x|a)-1)-\mu F_a(x|a)}{f(x|a)} & \text{if} \  x>\hat{x}.\\
\end{cases}
\label{e12}
\end{equation}
For simplicity, we assume for the next result that  $f_{Y|A}(y|a)= \mathcal N(m(a),{\sigma_y} ^2)$ and $f_N(n)=\mathcal{N}(0,\sigma_{n}^{2}).$ It is worth pointing out that assumptions $\mathcal A_{1}-\mathcal A_{4}$ hold in this case (for Gaussian distribution). By definition, $x=y-n$ will be a Gaussian random variable  and $f_{X|A}(x|a)= \mathcal N(m(a),{\sigma_x} ^2)$. Notice when there is no error $x=y$ and $f_{X|A}(z|a)=f_{Y|A}(z|a)=\mathcal N(m(a),{\sigma_y} ^2)$, in the case of noise; however, $x=y-n$ so $f_{X|A}(z|a)=\mathcal N(m(a),{\sigma_x} ^2={\sigma_y} ^2-\sigma_{n}^{2})$. Accordingly, suppose  $f_{Y|A}(z|a)$ and $f_{X|A}(z|a)$ represent the pdf of $x$ in the absence and presence of noise.
Comparing $f_{X|A}(z|a)$ and $f_{Y|A}(z|a)$ for a variable $0\leq z\leq 1$,  we obtain:
\begin{equation}
\begin{cases}
f_{X|A}(z|a)\geq f_{Y|A}(z|a)&\text{if}\ \vert z\vert\leq c\\
f_{X|A}(z|a)<f_{Y|A}(z|a) &\text{if}\ \vert z\vert > c,\\
\end{cases}
\label{e7}
\end{equation}

where $c=m(a)+\sigma_x \sigma_y\sqrt{\frac {2\ln(\frac{\sigma_x}{\sigma_y})}{\sigma_x^2-\sigma_y^2}}$. 
\begin{cor}
\label{pou}
Suppose the profit $x$ is estimated with an estimation error. If $c=m(a)+\sigma_x \sigma_y\sqrt{\frac {2\ln(\frac{\sigma_x}{\sigma_y})}{\sigma_x^2-\sigma_y^2}}$ and the constraint multiplier $\mu$ in \eqref{e6} is $0$, the optimal contract induces the customer to do less underreporting (in the sense that the customer under-reports for a narrower range of load reduction) in the presence of estimation error as compared  to the case without error.
\end{cor}
\begin{proof}
Given the distribution assumptions on $y$ and $n$, $x=y-n$ will be a Gaussian random variable and its variance will be less than ${\sigma_y} ^2$. Therefore, $F_{Y|A}(z|a)<F_{X|A}(z|a)$ and $\vert F_{Y|A}(z|a)-1\vert>\vert F_{X|A}(z|a)-1\vert$. Based on \eqref{e7}, it can be noted that if $c>1$,  for $0\leq z\leq1$, $f_{Y|A}(z|a)< f_{X|A}(z|a)$. Thus, if $\mu=0$, 
\begin{equation}
\frac{\vert F_{Y|A}(z|a)-1\vert}{f_{Y|A}(z|a)}>\frac{\vert F_{X|A}(z|a)-1\vert}{f_{X|A}(z|a)}, 
\end{equation}
Therefore, comparing \eqref{e12}  and \eqref{e6} for $\mu=0$ and $\lambda<1$, the customer does less underreporting when there exists noise in the estimation in an optimal contract.
\end{proof}
\begin{remark}
 Comparing the two cases, we see that $\hat{x}$ is identical in the two cases. However, $x^+$  will decrease in the case when $\sigma_{n}^{2}>0$. .
 \end{remark}

\subsection{Optimal Share Allocated to the Customer}
The following result shows that the optimal contract must utilize the option of giving shares to the customer.
\begin{cor} The value of $\alpha$ is strictly positive  in the optimal contract.
\end{cor}
\begin{proof}
Differentiating $\Pi$ with respect to $\alpha$ yields
\begin{equation}
\Pi'(\alpha)=E(R-x)-\int_{\hat{x}}^{1}\lambda(1+R_\alpha)(F-1) dx,
\end{equation} which can be reduced to
\begin{equation}
\Pi'(\alpha)=E(R-x)
\label{sl}
\end{equation}
This implies that  if the expectation of the distortion of the load  is positive (respectively negative), $\Pi$ will be increasing (respectively decreasing) with respect to $\alpha$. If $\alpha=0$, $R-x$ is equal to $0$ for $x\leq\hat{x}$ and positive for $x>\hat{x}$ (based on assumption $\mathcal{A}_{4}$ and its generalization to $y$ in Lemma 1, 
\begin{equation}
\frac{(1-\lambda)(F_{Y|A}(x|a)*f_{N}(x)-1)-\mu{(F_{Y|A})}_a(x|a)*f_{N}(x)}{f_{Y|A}(x|a)*f_{N}(x)}
\end{equation}
 is strictly concave). Thus, given that $R-x$ is continuous, $E(R-x)$ is strictly positive. As $\alpha$ increases,  \eqref{e6} indicates that the curve of $R-x$ shifts down, so that $R-x=-\alpha$ for $x<\hat{x}$. Consequently, $E(R-x)$ decreases as $\alpha$ increases. Thus, for a large enough $\alpha$, we have that $E(R-x)=0$. For this critical value of $\alpha$, \eqref{sl} implies that $\Pi'(\alpha)=0$. Further, this is clearly a maxima.
\end{proof}

\section{CONCLUSIONS}
In this paper,  we designed an optimal contract between a demand response aggregator (DRA) and a customer for incentive-based demand response. In this set up, the DRA asks  the customer to reduce her demand and compensates her for this reduction. However, since the DRA must supply the customer with as much power as she desires, a strategic customer can temporarily increase her base load to report a larger reduction after the demand response event. Based on management contract design problem, we proposed an optimal contract  that maximizes DRA's utility by incentivizing the customer both to make the maximal effort in reducing the load and  not to falsify the base load. The proposed optimal contract consists of two parts: a share of the DRA's profit in demand response event and  a part that is compensation paid to customer depending on load reduction as measured. Further, some properties of the customer share of the profit and the compensation paid to her were discussed.

Future work will involve considering the dynamic problem, impact of pricing, and also the multiple  customers and ownership case. Relating this work to the game theoretic set ups in~\cite{imp7} and~\cite{imp9} is also of interest.

\addtolength{\textheight}{-12cm}   




\section*{ACKNOWLEDGMENT}

We  would like to thank Dr. Thomas A.Gresik from Department of Economics in the University of Notre Dame for his insights and comments.


\bibliographystyle{IEEEtran}
%

\bibliography{a}

\enlargethispage{3\baselineskip}

\end{document}